\documentclass[a4paper]{article}

\usepackage{graphicx}
\usepackage{latexsym}
\usepackage{amsmath}
\usepackage{amssymb}

\makeatletter
\renewcommand{\section}{\@startsection
  {section}{1}{0mm}{-1.3\baselineskip}{0.7\baselineskip}%
  {\normalsize\bfseries}}
\makeatother

\newcommand{\paren}[1]{\ensuremath{\left(#1\right)}}

\newcommand{\e}{\mathrm{e}}

\newcommand{\pswitch}{P_{\mathop{\mathrm{switch}}}}

\begin{document}

\title{Prize insights in probability, and one goat of a recycled error\\\vspace{1ex}%
\normalsize Jason Rosenhouse.  \emph{The Monty Hall Problem: the remarkable story of math's most contentious brain teaser.}  Oxford University Press, 2009.  RRP US\$24.95.}
\author{Anthony B. Morton\\\small Honorary Fellow, Melbourne School of Engineering, University of Melbourne}
\date{November 2010}
\maketitle

\noindent
Jason Rosenhouse's `great big Monty Hall book' \cite{r:tmhp} is published by Oxford; but it was on one of my occasional pilgrimages to Cambridge that I stumbled upon it.
As one who has taken an unnatural interest in this subject over the years, I found it unsurprising that enough material existed out there to fill a book on the Monty Hall problem, but Rosenhouse is to be commended for having taken the trouble to write it.
And he has done so with suitable breadth, from the raw maths of the problem and others like it, to its psychological and philosophical implications.

The Monty Hall problem is the TV game scenario where you, the contestant, are presented with three doors, with a car hidden behind one and goats hidden behind the other two.
After you select a door, the host (named Monty Hall after the actual host of \emph{Let's Make a Deal}, a 1960s American TV show) opens a second door to reveal a goat.
You are then invited to stay with your original choice of door, or to switch to the remaining unopened door, and claim whatever you find behind it.
Assuming your objective is to win the car, is your best strategy to stay or switch, or does it not matter?
Rosenhouse provides the definitive analysis of this game, along with several intriguing variations, and discusses many of the wider issues that emerge.

Some reviewers, such as David Spiegalhalter in the \emph{LMS Newsletter} \cite{ds:tmhp} and Donald Granberg in \emph{Science} \cite{g:tdgtmhp} have queried the mathematical level of the book.
Indeed the maths content is in a seldom-acknowledged middle place: uncomfortably informal for those used to an academic standard of rigour, but still undeniably heavy for those expecting a popular treatment.
Yet the style of presentation is not as off-beat as might be thought by either `experts' or non-mathematicians: it will be quite familiar to engineers, for example.
The ``contention'' in Rosenhouse's title clearly doesn't reside on one side or other of C.P.\ Snow's Two Cultures---where one is either a robotic theorem-prover or a learned hater of equations---and nor can a book like this.
Those who take enough interest in brainbenders like Monty Hall to read books about it tend to be generalists who've kept closely in touch with their high school or undergrad mathematics.
Recreational maths tragics like myself, in other words.

I've certainly had fun digesting \emph{The Monty Hall Problem} and engaging critically with its contents.
The following---a much longer article than I had intended, for which I apologise in advance---is an attempt to share some of my reflections on the book's contents, as well as to detail just one or two criticisms.

\section*{Probability as Logic}

Part of what drives people to read or to write books about puzzles like Monty Hall is that there is a serious side beneath the sheer entertainment value of arguing about strategy in TV game shows.
To a diehard Bayesian like myself, whose take on probability theory owes much to having read Edwin Jaynes' posthumous text on the subject \cite{j:ptls}, the Monty Hall problem is a made-to-order demonstration of probability theory as the logical analysis of information---or `epistemic probability', as Rosenhouse later describes it in his Chapter 7.

As Rosenhouse makes very clear, the problem is fundamentally all about information.
In most versions of the problem, Monty is presumed to know what is behind each door, making clear that the only `uncertainty' involved is your own lack of information about where the car is.
You will assign initial probabilities based on this lack of information, which should then be updated as information is revealed over the course of the game.
On the `logical Bayesian' view expounded by the Jaynes school, the probabilities do not just represent your own subjective beliefs: they are definite quantities that measure the plausibility of a proposition such as ``the car is behind door 1'' on the basis of specific information.

Consider the situation at the outset.
A fundamental assumption is that you are completely ignorant as to where Monty has concealed the car, except that it is behind one of the three doors.
What this means more precisely is that if (in some parallel universe) you were presented with the same three doors but in a different order, this would not affect your probability assignments.
In other words, your background knowledge gives you no reason to rate any one door as more likely than any other to conceal the car: the doors are `exchangeable'.
The only rational course is then to assign the same probability to each door, and since there are three doors and you know the car is certainly behind one of them, you give each a probability of $1/3$.

As Rosenhouse explains fully in his Chapter 3, the correct procedure for updating these probabilities in the light of new information, such as revealed by Monty opening a door, is given by rules including Bayes' Theorem and the Law of Total Probability.
Rosenhouse presents these rules in the way one is accustomed to seeing in high school: conditional probability is introduced as a subordinate concept defined in terms of joint and single probabilities, and it is shown how the formulae that result are intuitively plausible.

Yet the reader, like many thoughtful high school students, might harbour some lingering doubts about the universal validity of such a postulated definition.
A neat resolution to this question was actually provided some decades ago, though it is not widely known.
First, note that it can be shown (the working is quite similar to that found in Rosenhouse's Chapter 3 and most high school textbooks) that the rules of basic probability all follow from two fundamental formulae, which we might call the `sum rule' and `product rule':
\begin{gather}
P(A | I) + P(\bar{A} | I) = 1,
\label{eq:sumrule} \\
P(AB | I) = P(A | I) P(B | AI).
\label{eq:prodrule}
\end{gather}
(Here I have followed the Bayesian convention of making all probabilities conditional on the term $I$, which represents a collection of propositions defining our `background information'.
$I$ may include, for example, the specification of the `probability space' which gives probabilities according to the standard Borel-Kolmogorov formulation of probability theory.
This inclusion of $I$ may seem just a needless distraction, but Rosenhouse will himself acknowledge the importance of this, right at the end of his book.
The notation $\bar{A}$ is used for the Boolean negation of $A$, or its set-theoretic complement when $A$ is a subset of a sample space.)

The physicist R.T. Cox showed in a 1961 monograph \cite{c:lpi} that the rules (\ref{eq:sumrule}) and (\ref{eq:prodrule}) emerge as theorems in a certain kind of multi-valued logic.
In this logic, any two propositions $A$ and $B$ determine a real number $P(A | B)$ which is taken to measure the degree of plausibility of proposition $A$ in light of the truth of proposition $B$.
(This of course corresponds to what we generally mean by conditional probability, but where $A$ and $B$ may be \emph{any} propositions we care to state.)
If the truth of $B$ implies that $A$ is in fact true (respectively, false), then $P(A | B)$ is given its largest (respectively, smallest) allowable value; these can be taken as 1 and 0 respectively with no loss of generality.
Crucially, in this special case of true and false propositions, Cox required his multi-valued logic to conform with ordinary Boolean logic.
(So for example, $P(AB | C)$ should be 1 whenever both $P(A | C) = 1$ and $P(B | C) = 1$, but should be zero when either of the latter are zero.)
It was also posulated on fairly reasonable grounds that, for example, if the truth of $C$ increases the probability of $A$ but leaves the probability of $B$ unchanged, then it should not reduce the probability of the joint proposition $AB$; but on the other hand it should reduce the probability of $\bar{A}$.
It turns out that these basic considerations suffice to prove (\ref{eq:sumrule}) and (\ref{eq:prodrule}), so that it is no longer necessary to postulate the latter \emph{a priori}.

Rosenhouse introduces his basic probability concept in Chapter 2.
Here, the probability $P(A)$ is defined for ``outcomes'' $A$ which are considered to be subsets of a predefined sample space.
A Bayesian would write this as $P(E_A | I)$, where $I$ states the definition of the sample space and the `experiment' being conducted, and $E_A$ is the proposition ``the outcome of the experiment belongs to the subset $A$ of the sample space.''
But the distinction is largely just one of syntax, and there is (mostly!) no difficulty in regarding $P(A)$ as a convenient shorthand and interpreting all of Rosenhouse's probabilities according to a strict Bayesian framework.

So why bother with the Bayesian framework at all, rather than just being content to define probabilities for subsets of an underlying sample space?
Most of the time, I would agree it doesn't matter: the notion of drawing from a sample space works admirably, for the most part, throughout the mathematical discussions in \emph{The Monty Hall Problem}.

Occasionally, however, it pays to have the extra flexibility to discuss arbitrary propositions.
I will consider one example in some detail.
Throughout all discussions of the `classical' Monty Hall problem, Rosenhouse gives great emphasis to the need for Monty to \emph{choose at random} between the two remaining doors in the case where he knows they both conceal goats.
This assumption is necessary to ensure that when we assign probability $1/2$ to each outcome in this situation, those probabilities do in fact correctly measure the way the sample space is carved up by Monty over many repetitions of the game.

But as I said at the outset, the Monty Hall problem is all about information.
Rosenhouse points out that most statements of the classical Monty Hall problem fail to include any explicit assumption about Monty's behaviour when he has two goats to choose from, merely stating that Monty will always open a door that conceals a goat.
Yet those who state the problem this way still manage to reason, mostly correctly, to the conclusion that switching wins the car with probability $2/3$.
And in any case, what does it mean to say that Monty ``chooses at random'' when we are interested in the outcome of just a single instance of the game?
Rosenhouse touches on this question in Chapter 7, when discussing some philosophical arguments about different interpretations of probability.
I'll have more to say about those arguments in due course.

For now, let us consider Proposition $X$: ``Monty knows that two doors both conceal goats, and will open one of them, but you do not know how he decides which to open.''
Proposition $X$ does not talk about subsets of a sample space, or about long-run frequencies in a random experiment, but it is a valid proposition of Cox logic.
So if we can find a reasonable way to assign unique numbers $P(A | X)$ and $P(B | X)$, where $A$ and $B$ postulate that Monty opens the first or second door respectively, then ordinary probability theory will still produce a reasonable and unique answer to our problem.

The answer that immediately suggests itself is the same one indicated in our initial situation, where we have three doors and no knowledge of which one conceals the car.
As in that case, we can argue that since we have no reason to believe Monty is more likely to choose one door than any other, the two doors are exchangable and hence we should put $P(A | X) = P(B | X) = 1/2$.

This seems quite reasonable, but perhaps not yet quite convincing.
Jaynes, in a key 1968 paper on prior probabilities \cite{j:pp}, provides an ingenious combinatorial argument that assigning equal probabilities really is the most appropriate and `conservative' choice in light of assumptions like Proposition $X$.
Consider a large number $N$ of identical trials in which Monty chooses either the first or the second door.
Suppose he opens the first door in $M$ of these trials, and the second door in the remaining $N - M$ trials.
Then the ratio $p = M/N$ is the statistical frequency which ought to closely match our probability assignment $P(A | X)$.

Jaynes' argument focusses on the \emph{number of ways} in which Monty can realise this frequency as the outcome of the $N$ trials.
For a given value of $M$, this is just the binomial coefficient $\binom{N}{M}$, out of a total of $2^N$ possible ways of choosing one or other door over the course of $N$ trials.
If $N$ and $M$ are sufficiently large, it turns out we can approximate this huge number of combinations as
\begin{equation}
\binom{N}{M} \approx \e^{N \cdot H(p)}
\label{eq:binapprox}
\end{equation}
where
\begin{equation}
H(p) = - p \log(p) - (1 - p) \log(1 - p), \qquad p = \frac{M}{N}.
\label{eq:entropy}
\end{equation}
The formula (\ref{eq:entropy}) is the \emph{Shannon entropy} of an event having probability $p$ (which is discussed by Rosenhouse in a slightly different form in Chapter 4).
$H(p)$ takes its maximum value when $p = 1/2$, which means that according to (\ref{eq:binapprox}), the frequency $p = 1/2$ is the one that can be realised in the greatest possible number of ways.
But what may not be apparent from (\ref{eq:binapprox}) is just \emph{how} great this greatest number is.
Suppose we have $N = 1000$ trials; then there are some $5.4 \times 10^{130}$ ways to achieve $p = 0.5$, but just $6.5 \times 10^{126}$ ways to achieve $p = 0.4$.
That is to say, there are some 10,000 more combinations giving $p = 0.5$ than there are giving $p = 0.4$.

It follows that if one is to assign probabilities $P(A | X)$ and $P(B | X)$ that express `maximum ignorance' about Monty's behaviour, there is a further compelling reason to set $P(A | X) = P(B | X) = 1/2$: this turns out to be consistent with the greatest possible range of future actions by Monty.
This `maximum entropy' rule readily generalises, and has proved to be a valuable tool in Bayesian data analysis.

\section*{Variations and Pitfalls}

Naturally, the above pertains to the `classical' Monty Hall problem, where Monty always contrives to open a door that conceals a goat, and your probability of winning by switching is $2/3$.
Rosenhouse devotes Chapter 3 to explaining why this probability actually drops to $1/2$ in the subtly different case where Monty selects a door to open at random, and happens to reveal a goat.
In this situation, unlike the classical game, there is a very real chance that Monty will reveal the car and end the game; the fact that he does not reveals information that is relevant to the location of the car.
We are now in the domain of elementary `Bernoulli sampling without replacement', where the probabilities depend only on the number and type of outcomes that remain available, independently of what has been drawn previously and of which door we choose.
It's not out of the question that because most of us learn about probability this way, this leads people to think the same should apply in the classical Monty Hall problem as well (and in Chapter 6, we see some experimental evidence along these lines).

In Chapters 4 and 5, Rosenhouse further develops the problem by analysing the manifold variations that have appeared in the mathematical literature.
Among these variations are some, originally discussed in a paper by Georges and Craine \cite{gc:gmd}, that consider the presence of multiple cars, and more generally, a variety of prizes of differing value, certain of which may or may not be revealed by Monty.
It is here that I suggest, with great trepidation, that Rosenhouse has committed his one serious mathematical error---or perhaps, recycled an error committed earlier by others.

Paulo Ventura Ara\'{u}jo, reviewing the book for the \emph{Newsletter of the European Mathematical Society} \cite{a:tmhp}, states that these mathematical variations ``are not the stuff of bedtime reading.''
Perhaps he did not have nutty engineers like me in mind.
In any case, it was while reading Chapter 5 in bed one night that I came across the following passage in Section 5.4 (emphasis mine):
\begin{quote}
This time we still have $n$ doors, but now there are $1 \leq j \leq n - 2$ cars and $n - j$ goats.
After making your initial choice, Monty opens one of the other doors \emph{at random}.
Should you switch?

Since there are $j$ cars and $n$ doors, we see that our initial choice conceals a car with probability $j / n$.
This probability \emph{will not change regardless of what Monty reveals}, and therefore represents the probability of obtaining a car by sticking.
\end{quote}
This was enough to jolt me awake.
After all, it had just been carefully argued in Chapter 3 that when Monty chooses a door randomly, this \emph{does} affect the posterior probability that our initial choice concealed the car!

Sure enough, we can repeat the calculation using Bayes' Theorem to show that the conclusion of Section 5.4 is erroneous: the probability of winning a car by switching to an unopened door does \emph{not} go up or down relative to the probability of winning by sticking, depending whether Monty's random choice reveals a car or a goat.
Instead, the effect of Monty's random choice is to redistribute the probabilities evenly over all the unopened doors---just as Chapter 3 concluded in the three-door, single-car case.
So the probability of winning a car in this case is $j / (n - 1)$ if Monty reveals a goat and $(j - 1) / (n - 1)$ if Monty reveals a car, no matter whether you stick with the original door or switch to any other door.

Here is the calculation.
Following Rosenhouse's notation, let $F_c$ and $F_g$ denote the propositions that your first choice conceals a car, respectively a goat, and $S_c$ and $S_g$ the propositions that you switch to a door concealing a car, respectively a goat.
Also, let $M_c$ and $M_g$ denote the propositions that Monty selects at random a door concealing a car, respectively a goat.
On our background information $I$ alone, before anything else has happened, we clearly have
\[P(F_c | I) = \frac{j}{n}, \qquad P(F_g | I) = \frac{n - j}{n}.\]
So far so good.
Now, suppose that Monty (choosing randomly) reveals a goat.
Since Monty does not open the door you have selected, the probabilities here depend on whether your first choice was good or bad.
In fact it is readily seen that
\[P(M_g | F_c I) = \frac{n - j}{n - 1}, \qquad P(M_g | F_g I) = \frac{n - j - 1}{n - 1}.\]
The conditional probabilities of winning a car by switching are correctly calculated by Rosenhouse as
\[P(S_c | F_c M_g I) = \frac{j - 1}{n - 2}, \qquad P(S_c | F_g M_g I) = \frac{j}{n - 2},\]
but now, to apply the Law of Total Probability we need to calculate $\pswitch$, the probability you win by switching, as
\[\pswitch = P(S_c | M_g I) = P(S_c | F_c M_g I) P(F_c | M_g I) + P(S_c | F_g M_g I) P(F_g | M_g I).\]
Thus we need the probabilities of your first choice concealing a car or a goat, conditional on the information that Monty revealed a goat.
These can be calculated using Bayes' Theorem as
\[\begin{split}
P(F_c | M_g I) &= \frac{P(F_c | I) P(M_g | F_c I)}{P(M_g | F_c I) P(F_c | I) + P(M_g | F_g I) P(F_g | I)} \\
   &= \frac{j}{n} \cdot \frac{n - j}{n - 1} \paren{\frac{n - j}{n - 1} \cdot \frac{j}{n}
      + \frac{n - j - 1}{n - 1} \cdot \frac{n - j}{n}}^{-1}
   = \frac{j}{n - 1}
\end{split}\]
and
\[P(F_g | M_g I) = \frac{P(F_g | I) P(M_g | F_g I)}{P(M_g | F_c I) P(F_c | I) + P(M_g | F_g I) P(F_g | I)}
   = \frac{n - 1 - j}{n - 1}.\]
Plugging these into the formula for $\pswitch$ gives
\[\pswitch = \frac{j - 1}{n - 2} \cdot \frac{j}{n - 1} + \frac{j}{n - 2} \cdot \frac{n - 1 - j}{n - 1} = \frac{j}{n - 1}.\]
As anticipated, both $\pswitch$ and $P(F_c | M_g I)$---the probability you win by sticking with your first choice---come out as $j / (n - 1)$, which is the number of remaining cars divided by the number of remaining doors.

The working is similar if Monty reveals a car.
Now we have
\[P(M_c | F_c I) = \frac{j - 1}{n - 1}, \qquad P(M_c | F_g I) = \frac{j}{n - 1},\]
and
\[P(S_c | F_c M_c I) = \frac{j - 2}{n - 2}, \qquad P(S_c | F_g M_c I) = \frac{j - 1}{n - 2}.\]
Again, we need the posterior probabilities of your first choice, conditional on $M_c$:
\[\begin{split}
P(F_c | M_c I) &= \frac{P(F_c | I) P(M_c | F_c I)}{P(M_c | F_c I) P(F_c | I) + P(M_c | F_g I) P(F_g | I)} \\
   &= \frac{j}{n} \cdot \frac{j - 1}{n - 1} \paren{\frac{j - 1}{n - 1} \cdot \frac{j}{n}
      + \frac{j}{n - 1} \cdot \frac{n - j}{n}}^{-1}
   = \frac{j - 1}{n - 1}
\end{split}\]
and
\[P(F_g | M_c I) = \frac{P(F_g | I) P(M_c | F_g I)}{P(M_c | F_c I) P(F_c | I) + P(M_c | F_g I) P(F_g | I)}
   = \frac{n - j}{n - 1}.\]
Finally, we have
\[\begin{split}
\pswitch &= P(S_c | M_c I) = P(S_c | F_c M_c I) P(F_c | M_c I) + P(S_c | F_g M_c I) P(F_g | M_c I) \\
   &= \frac{j - 2}{n - 2} \cdot \frac{j - 1}{n - 1} + \frac{j - 1}{n - 2} \cdot \frac{n - j}{n - 1} = \frac{j - 1}{n - 1}.
\end{split}\]
Once again, the probability of winning a car is equal to the number of remaining cars divided by the number of remaining doors---here $(j - 1) / (n - 1)$---regardless of whether you switch or not.

Of course, this raises the question: what conditions are necessary in this more general case to ensure your probability of picking a car the first time \emph{is} unchanged after Monty opens another door (and thereby make the original reasoning in \cite{gc:gmd} valid)?
The answer: only if the probability Monty reveals a goat \emph{does not change} according to whether your first chosen door conceals a car or a goat.
To see this, we use Bayes' Theorem again.
The condition we require is that
\[P(F_c | M_g I) = P(F_c | I).\]
From our formula for $P(F_c | M_g I)$ above, this condition is equivalent to saying that
\[P(M_g | F_c I) = P(M_g | F_c I) P(F_c | I) + P(M_g | F_g I) P(F_g | I),\]
or
\[P(M_g | F_c I) \paren{1 - P(F_c | I)} = P(M_g | F_g I) P(F_g | I).\]
But since your first choice of door is certain to conceal either a car or a goat, we surely have
\[1 - P(F_c | I) = P(F_g | I).\]
Substituting this above and simplifying (assuming that $P(F_g | I) > 0$) leads immediately to
\[P(M_g | F_c I) = P(M_g | F_g I).\]
Note that if there is just one car (which it is possible for you to choose at the outset) then this condition effectively forces Monty to always reveal a goat.
In general, if Monty behaves so as always to reveal a goat, this suffices to ensure your probability of having picked the car does not change when Monty reveals it.
But if there are multiple cars, it is possible for Monty to behave in other ways---including to always reveal a car---and still leave your probability unchanged.
A close reading of \cite{gc:gmd} suggests that Georges and Craine were assuming Monty to be following one or other of these predetermined strategies, though this was not made altogether explicit.
But if Monty simply chooses a door at random, then the above condition will always be violated, since random Monty is always slightly more likely to reveal a goat if your chosen door conceals the car.

(Interestingly, the same paper by Georges and Craine, with its incompletely specified problem, also features in a recent e-book by another of my favourite semi-popular mathematical writers, Julian Havil \cite{h:isscc}.
Havil's book is concerned with much more wide-ranging fare, but Chapter 6 does include a discussion of the Monty Hall problem, including the same Georges--Craine inspired variation and the same problematic conclusion: that Monty increases your chance of winning by switching by revealing a goat and reduces it by revealing a car, which is incorrect when he chooses at random.)

The good news is that the next piece of analysis in Section 5.4, where Monty is assumed to reveal a goat with probability $p$ and a car with probability $1 - p$, all goes through correctly---precisely because the probability $P(M_g | I)$ is \emph{assumed} equal to $p$ irrespective of your first choice of door.
Therefore the above condition is fulfilled, and the prior and posterior probabilities of $F_c$ and $F_g$ are identical.
Likewise, in the multiple-door scenario that follows, it is assumed that Monty will always open $m$ doors to reveal $k$ cars, thus Monty's behaviour is the same irrespective of whether $F_c$ or $F_g$ is true and the analysis goes through as Rosenhouse describes.

In Section 5.6, Rosenhouse moves on to consider the most general Monty Hall scenario (again drawing on \cite{gc:gmd}), with an arbitrary number of doors and an arbitrary number of prizes with arbitrary values.
Here it appears we may have run into trouble once again, because in this variation Monty is once again assumed to choose doors at random, yet Rosenhouse states that ``these probabilities [of winning a specific prize with our first guess] do not change when Monty opens a door.''

Because the prizes now have differing values, the winnings are now framed in terms of expected values over the long run.
More precisely, Rosenhouse is dealing with \emph{conditional expectations}, which obey their own form of the Law of Total Probability.
Unfortunately, Rosenhouse (like the authors in \cite{gc:gmd}) calculates with these values without explicitly stating they are conditional, or what the exact conditioning proposition is.
In fact, all expectations calculated for both the `stick' and `switch' strategy are conditional on the type of prize revealed by Monty; accordingly, the weighting probabilities need to be conditional on what Monty reveals as well.
Only if Monty decides \emph{in advance} which prize to reveal, irrespective of your choice, does the analysis go through as stated in \cite{gc:gmd} and Section 5.6.

For the `stick' strategy, it is certainly true that if each door conceals a prize (where some `prizes' may have zero value) and the mean value of all prizes is $V$, then your expected win from sticking is $V$, averaged over \emph{all} possible plays of Monty's game (and irrespective of what Monty reveals or how Monty behaves in each instance).
However, when Monty chooses at random, this (unconditional) expected value is itself the average of very \emph{different} conditional expectations, that depend on what Monty actually reveals.

Consider again the scenario with three doors, one car (of value $3V$) and two goats (of value zero).
If Monty's random selection reveals a goat, then you win a car half the time \emph{out of those specific games} where Monty revealed a goat, and your expected winnings \emph{from those specific games} are $3V/2$.
But random Monty will only reveal a goat $2/3$ of the time; in the other $1/3$ of cases Monty reveals the car and your winnings are zero.
Sure enough, when you average out the expected winnings of $3V/2$ two-thirds of the time and zero one-third of the time, the overall expected winnings are $V$---which is the same as your expected winnings from sticking in the classical scenario, where Monty always reveals a goat, and you win the car $1/3$ of the time.

The situation is similar when calculating expected values for the switching strategy.
Rosenhouse's concrete example is a six-door game where two doors conceal cars each worth \$20,000, two doors conceal motorcycles each worth \$10,000, and two doors conceal refrigerators each worth \$300.
Let $F_c$, $F_m$ and $F_r$ denote the propositions that your first choice of door conceals a car, a motorcycle and a refrigerator respectively, and let $M_c$ be the proposition that Monty opens a door at random and reveals a car.
Rosenhouse calculates the conditional expectations for the value $V_S$ of the second choice as
\[E[V_S | F_c M_c I] = \text{\$5,150}, \qquad E[V_S | F_m M_c I] = \text{\$7,650}, \qquad
E[V_S | F_r M_c I] = \text{\$10,075},\]
and then determines the overall expectation conditional on $M_c$ as
\[E[V_S | M_c I] = \frac{1}{3} (\text{\$5,150}) + \frac{1}{3} (\text{\$7,650}) + \frac{1}{3} (\text{\$10,075})
= \text{\$7,625}.\]
It is in the latter calculation that a mistake has occurred, because (if Monty really has chosen at random) the weights ought to be the posterior probablities like $P(F_c | M_c I)$, but instead the calculation has used the prior probabilities like $P(F_c | I)$.
We can correct this easily using Bayes' Theorem to calculate the posterior probabilities as
\[P(F_c | M_c I) = \frac{1}{5}, \qquad P(F_m | M_c I) = \frac{2}{5}, \qquad P(F_r | M_c I) = \frac{2}{5},\]
in which case the correct expected value is
\[E[V_S | M_c I] = \frac{1}{5} (\text{\$5,150}) + \frac{2}{5} (\text{\$7,650}) + \frac{2}{5} (\text{\$10,075})
= \text{\$8,120}.\]

In the general case, there are assumed to be $m$ doors and $n_i$ prizes each worth $v_i$, where $\sum_i n_i = m$.
We set $t = \sum_i (n_i v_i)$; this is the same as $mV$ where $V$ is the mean value of all prizes.
The expected value of the prize when switching, conditional on $F_i$ (``your first choice of door conceals a prize of value $v_i$'') and on $M_r$ (``Monty's randomly revealed prize has value $v_r$'') is correctly calculated as
\[E[V_S | F_i M_r I] = \frac{t - v_i - v_r}{m - 2}.\]
Now we require the posterior probabilities $P(F_i | M_r I)$.
Each of these is given by Bayes' Theorem as
\[P(F_i | M_r I) = \frac{P(F_i | I) P(M_r | F_i I)}{\sum_k P(F_k | I) P(M_r | F_k I)}.\]
For all $i$ we have $P(F_i | I) = n_i / m$, while the value of $P(M_r | F_i I)$ is $n_r / m$ if $i \neq r$ and $(n_r - 1) / m$ if $i = r$.
It follows that
\[\sum_k P(F_k | I) P(M_r | F_k I) = \sum_k \frac{n_k}{m} \frac{n_r}{m} - \frac{n_r}{m^2}
   = \frac{n_r}{m} \frac{m - 1}{m}\]
and so
\[P(F_i | M_r I) = \begin{cases} \frac{n_i}{m - 1} & \text{if $i \neq r$}, \\
   \frac{n_i - 1}{m - 1} & \text{if $i = r$}. \end{cases}\]
Finally, we have
\[\begin{split}
E[V_S | M_r I] &= \sum_{i} \frac{n_i}{m - 1} \frac{t - v_i - v_r}{m - 2} - \frac{1}{m - 1} \frac{t - 2v_r}{m - 2} \\
   &= \frac{1}{(m - 1) (m - 2)} \paren{\sum_{i} n_i (t - v_r) - \sum_i n_i v_i + (2v_r - t)} \\
   &= \frac{(m - 2) (t - v_r)}{(m - 1) (m - 2)}
   = \frac{t - v_r}{m - 1}.
\end{split}\]
So much for the expected value of switching: what is it now for the `stick' strategy?
It is not simply $V$, as Rosenhouse assumes, because the possible values $v_i$ for our first choice $F_i$ must now be weighted by the posterior probabilities $P(F_i | M_r I)$.
Let $E[V_F | M_r I]$ denote the expected value of our first choice, given that Monty randomly reveals a prize of value $v_r$: then
\[E[V_F | M_r I] = \sum_{i} \frac{n_i}{m - 1} v_i - \frac{1}{m - 1} v_r = \frac{t - v_r}{m - 1} = E[V_S | M_r I].\]
In other words, the conclusion is just as found in the analysis of Section 5.4 above: Monty's random choice has the effect of redistributing the probabilities equally over all unopened doors.
With Monty choosing at random, it doesn't matter whether you switch or stick: your expected winnings are the same in each case.
This is Bernoulli sampling without replacement, just as before.

Notice too that the expected values $E[V_F | M_r I] = E[V_S | M_r I]$ become equal to $V = t / m$ precisely when $v_r = V$: exactly the same as found by Rosenhouse and the authors of \cite{gc:gmd} when Monty uses a predetermined strategy.
It follows that if the revealed prize has value $v_r < V$, the expected values for both the `stick' and `switch' strategy become greater than $V$, while if $v_r > V$, both expected values become less than $V$.
But one can see intuitively why this should be so: once random Monty reveals and eliminates a prize, it is as if you are facing an entirely new game based on the values of the remaining prizes.
It is very different to the situation where Monty is presumed to always select a prize of value $v_r$ regardless of your initial choice.
Then, it is only the expected value $E[V_S | M_r I]$ that changes, as Rosenhouse and \cite{gc:gmd} describe: up when $v_r < V$ and down when $v_r > V$.

And now, a confession: in an earlier version of this review article, I gave the value of $E[V_S | M_r I]$---the expected value of switching---as above, but failed to notice that Monty's random selection will also change the posterior value of $E[V_F | M_r I]$, the expected value for sticking.
I too was distracted by the clear fact that $E[V_F | I]$, the expectation over all games, is equal to $V$, temporarily forgetting that the values sought are the ones conditional on $M_r$.
It therefore appeared that the original conclusion---switch when $v_r < V$ and stick when $v_r > V$---would still hold in the case where Monty opens a door at random.
Alas, it is not so, and it took me several days to spot my own error.

This episode starkly reveals the importance of stating one's assumptions very clearly when discussing the multitude of ways in which Monty's game can be varied.
Unlike many other authors, Rosenhouse is quite clear on this point.
But the many hidden subtleties in the Monty Hall scenario are such as to confound even the most vigilant---myself included!

\section*{Philosophical Entanglements}

Chapter 7 of Rosenhouse's book rounds out the discussion by considering some philosophical conundrums around the distinction that is often drawn between so-called `epistemic' and `statistical' probability.
It is central to the Bayesian view that outside the realm of thought experiments, and leaving aside debates about the interpretation of quantum mechanics, all `real world' probabilities are `epistemic' in the sense Rosenhouse means here.
Furthermore, as long as one is always careful to update one's `epistemic' probability assignments as any new information comes to light, these probabilities in the long run should always coincide with `statistical' observed frequencies in multiple runs of an experiment.
(Rosenhouse notes that many authors prefer to use the term `epistemic statistical probability' when the new information comes in the form of data from statistical trials, but the principle is a general one.)
If any discrepancy remains between probabilities and observed frequencies, this indicates to a Bayesian the presence of some relevant information that has not properly been taken into account.
Accordingly, any suggestion that epistemic and statistical probabilities should in principle differ should be closely scrutinised.

Rosenhouse reproduces an argument of Moser and Mulder \cite{mm:pirdm} which purports to find a difference between the appropriate `statistical' probability in long runs of the Monty Hall game, and the appropriate `epistemic' probability in a single instance of the game.
The argument hinges on our ignorance of Monty's procedure for revealing a goat; an issue that was extensively discussed earlier in this article.
The conclusion reached there is that faced with a lack of information, you should assign the probability that gives Monty maximum `freedom' to act, in the sense of our maximum entropy rule.
This assignment of probability leads to the $2/3$ probability of winning by switching, and this remains true whether we are considering a single game or a long run of similar games.

But this reasoning may appear to break down when Rosenhouse illustrates the Moser-Mulder argument with a slightly altered Monty Hall game.
In this variation the doors are numbered; and when Monty has to choose which door to open because both his options conceal goats, he acts in accordance with one of two hypotheses: the first (say $H_1$) is that Monty selects a door randomly; the second ($H_2$) is that Monty always selects the higher-numbered door.
There is also a \$100 bonus paid whenever you do not switch, to break the tie that occurs when you decide the car is equally likely to be behind the two unopened doors.

The crux of the matter is this: it's been shown in detail that on hypothesis $H_1$, you win the car $2/3$ of the time by switching, and your best strategy is always to switch.
On hypothesis $H_2$, you also win the car $2/3$ of the time by switching; but here you also have the opportunity to improve your winnings (slightly) by using what you know about Monty's behaviour.
If for example you initially pick door 3 and Monty opens door 1, then on $H_2$ this signals that Monty was forced to open door 1 because door 2 conceals a car (since Monty would not open door 1 otherwise), and you win the car with certainty by switching.
This or an equivalent scenario will occur $1/3$ of the time.
The other $2/3$ of the time, if you choose door 3 then Monty will open door 2, either because door 1 conceals a car or because Monty has a choice of doors and applies the `higher-numbered door' rule.
Each of these alternatives occurs with probability $1/2$ (out of those occasions where the car is not behind door 2), so by staying with door 3 you win the car half the time (just as you do by switching), but also pocket the \$100 bonus on each of these occasions.
Your best strategy on $H_2$, in other words, is: ``switch if Monty opens the lower-numbered of the two doors available to him; otherwise, stay with your initial choice.''

It follows that in some particular cases---say, where you initially select door 3 and Monty opens door 2---your best action is different depend on whether $H_1$ or $H_2$ is true.
This is despite $H_1$ and $H_2$ both giving the same `statistical' probabilty $2/3$ to winning by switching.
As Rosenhouse states:
\begin{quote}
Both versions lead to the same statistical probabilities for winning by switching.
But they mandate different behaviour in the individual case.
\end{quote}
All the same, this conclusion is not quite right, because it is not really true that the `statistics' are entirely the same for each hypothesis.
$H_1$ and $H_2$ do both lead to a uniform $2/3$ probability across all cases.
But $H_2$ additionally partitions the possible cases into two subsets, which you can tell apart as soon as Monty opens a door.
And in one of these subsets the probability for winning by switching is 1, while in the other it is $1/2$.
This is precisely what leads you to behave differently in different cases.
So Rosenhouse, quite correctly, dismisses the notion put by Moser and Mulder that there are somehow cases in which the rational course of action in a single instance actually differs from the rational course of action based on long-run considerations.

Also surviving the challenge is my own immodest contention above, that (epistemic) probabilities and observed frequencies will always converge on each other in the long run provided one accounts for all available information.
It is true enough that a statistical observed frequency, by itself, will not always yield sufficient information to identify the best strategy in any individual case; but this is quite different from saing that `epistemic' and `statistical' probabilities can fail to agree in the long run.
Notice after all that if Monty really were behaving according to $H_2$, then we could take the observed outcomes and also partition them into our two recognisable subsets.
And in each of these subsets, as well as the entire population of cases, the observed frequency will be found to match the probability given by Bayesian calculus.
Rosenhouse quotes a response by Horgan \cite{h:lmad} to Moser and Mulder, in which Horgan argues similarly for the equivalence of epistemic and `epistemic statistical' probabilities; I would go further, and contend that the asymptotic equivalence also extends to actual observed frequencies.

The other controversy featured in Chapter 7 has, in both Rosenhouse's view and my own, a much simpler resolution.
Baumann \cite{b:tdtpscp} describes a two-player variation on the Monty Hall problem, where two players select doors without each other's knowledge, and Monty then opens a door based on both their selections, before offering both the opportunity to switch.
Baumann's analysis indicates that for either player, there is a $3/7$ probability of winning a car by staying with their original choice and a $4/7$ probability of winning by switching.
(An alternative line of argument by Levy \cite{l:bomhpscp} is dismissed by Rosenhouse, correctly in my view.)

Baumann wishes to conclude from his correct analysis of the two-player problem that it leads to a paradox, on the basis that the two players are supposedly faced with identical information and yet can assign different probabilities to propositions such as ``the car is behind door 2''.
Here, Rosenhouse correctly identifies the fallacy in this reasoning: the two players are in fact not faced with identical information, because they have selected different doors to start with, and neither is provided with information on which door the other has selected.
Despite Baumann's objection, the part of each player's background information that stipulates their initial selected door is relevant to the players' probabilities, because it is logically connected to them via Monty's actions in response to the players' choices.
This, of course, underlines the importance of the background information $I$ in a probability symbol like $P(A | I)$: as Jaynes \cite{j:ptls} explains, many a `paradox' has been generated by failing to recognise that two people are generating probabilities from different prior information.

\vspace{0.5\baselineskip}
\hfill\hrulefill\hfill
\vspace{0.5\baselineskip}

\noindent
There ends Rosenhouse's tour of the Monty Hall problem and my own reflections thereon.
Far from its humble TV game show origins, the problem has acquired formidable power, so that even world-leading mathematicians such as Paul Erd\"{o}s have succumbed to it.
Together with its numerous cunning variations, the Monty Hall problem has astounding capacity to trip the wary and the unwary alike.
But study of the problem also proves to have its rewards, yielding precious insights into the subtle workings of probability theory.
And that's worth a lot more than a goat.

\bibliographystyle{plain}
\bibliography{bayes}

\end{document}